\input amstex
\documentstyle{amsppt}
\input epsf
\TagsOnRight
\pagewidth{5.5 in}
\pageheight{8 in}
\magnification=1200
\NoBlackBoxes
  
\def\e{{\roman e}}

\def\conj#1{\overline{#1}}
\def\inner#1#2#3{ \langle #1,\, #2 \rangle_{_{#3}}}
\def\Dinner#1#2{ \langle #1,\, #2 \rangle_{D}}
\def\Dpinner#1#2{ \langle #1,\, #2 \rangle_{D^\perp}}

\def\vartau{\tau\hskip2pt\!\!\!\!\iota\hskip1pt}

\def\FBa{1}
\def\FBb{2}
\def\BEEKa{3}
\def\BEEKb{4}
\def\BK{5}
\def\EE{6}
\def\EL{7}
\def\AK{8}
\def\SL{9}
\def\MRa{10}
\def\MRb{11}
\def\MRc{12}
\def\KR{13}
\def\SWa{14}
\def\SWb{15}
\def\SWc{16}
\def\SWd{17}
\def\SWe{18}
\def\WW{19}

%
%

\topmatter

\title
Fourier Invariant Partially approximating subalgebras of the irrational rotation C*-algebra
\endtitle

\rightheadtext{Partially approximating invariant subalgebras}

\author
S.~Walters
\endauthor

\affil
\sevenrm The University of Northern British Columbia
\endaffil

\date
\sevenrm May 30, 2001
\enddate

\address
Department of Mathematics and Computer Science, The University
of Northern British Columbia, Prince George, B.C.  V2N 4Z9  CANADA
\endaddress

\email
walters\@hilbert.unbc.ca \ \ or \ \ walters\@unbc.ca  \hfill
\break\indent{\it Website}: http://hilbert.unbc.ca/walters
\endemail

\thanks
Research partly supported by NSERC grant OGP0169928
\hfill {\sevenrm (\TeX File: pointproj.tex)}
\endthanks

\keywords
C*-algebras, irrational rotation algebras, automorphisms, inductive limits,
K-groups, AF-algebras, Theta functions
\endkeywords

\subjclass
46L80,\ 46L40,\ 46L35
\endsubjclass

\abstract 
For all transcendental parameters, the irrational rotation algebra is shown to 
contain infinitely many C*-subalgebras satisfying the following properties.  
Each subalgebra is isomorphic to a direct sum of two matrix algebras of the same 
(perfect square) dimension; 
the Fourier transform maps each summand onto the other;
the corresponding unit projection is approximately central;
the compressions of the canonical generators of the irrational rotation algebra 
are approximately contained in the subalgebra.
\endabstract

\endtopmatter


\document

\specialhead \S1. Introduction \endspecialhead

For each irrational number $\theta$, the irrational rotation C*-algebra $A_\theta$ is
generated by unitaries $U,V$ satisfying the Heisenberg commutation relation  $VU=\lambda
UV$, where $\lambda=e^{2\pi i\theta}$.   The Fourier transform on $A_\theta$ is the 
order four automorphism $\sigma$ (or its inverse) of $A_\theta$ satisfying $\sigma(U)=V$
and $\sigma(V)=U^{-1}$. It can be viewed as a non commutative analogue of the
Fourier transform of classical analysis (which itself plays an important role in the 
theory). Its square is the so-called flip automorphism and is fairly well understood 
(see [\FBa], [\BEEKa], [\BEEKb], [\BK], [\AK], [\SWa]).
An open problem of Elliott (private communication), is if the Fourier transform has an 
inductive limit structure.  More precisely, can $A_\theta$ be approximated by Fourier
invariant basic building blocks consisting of finite direct sums of circle algebras and 
points?  
If it can be shown that the answer is affirmative, then it would seem that one is only 
a small step away from proving that the fixed point subalgebra $A_\theta^\sigma$ is an 
AF-algebra.
One of the main difficulties in this problem (aside from some number theoretic ones), 
and which separates it from the corresponding flip problem, is that the projections 
needed in the construction are somewhat more exotic, as they can no longer be Rieffel 
projections.  Therefore, the author believes that a detailed study of such projections
is needed, and this was begun by Boca [\FBb] in his construction of Fourier invariant
projections, and later by the author in [\SWe] where projections orthogonal to their 
Fourier transform are constructed (such projections are needed in constructing direct 
sums of two identical circle algebras in which the Fourier transform swaps the two 
summands).  In [\SWd], a clue is given of the inductive limit structure of the Fourier 
transform by the construction of a (noncanonical) model of an inductive limit 
automorphism of order four on $A_\theta$ that agrees with $\sigma$ on $K_1$ 
(a pseudo-Fourier transform) such that the fixed point subalgebra is an AF-algebra.  
In [\SWb] and [\SWc], it was shown for the crossed product $A_\theta\rtimes_\sigma\Bbb Z_4$
that there is always an inclusion $\Bbb Z^9 \hookrightarrow K_0$, which, for a dense
$G_\delta$ set of $\theta$'s, is an isomorphism (and $K_1=0$).
The main thrust of this paper is to begin tackling Elliott's problem more directly 
(and canonically) by showing that one can construct infinitely many partially 
approximating finite dimensional subalgebras of $A_\theta$ (i.e., points) that are
invariant under the Fourier transform, as well as invariant direct sums of two circle 
subalgebras.  More precisely, the following result is proved. 


\proclaim{Main Theorem}
Let $\theta$ be a transcendental number.  Then the irrational rotation algebra 
$A_\theta$ contains infinitely many C*-subalgebras $\Cal M_q^2 = 
\Cal M_q\oplus \Cal M_q$ (orthogonal sum), where $\Cal M_q$ is isomorphic to 
$M_q(\Bbb C)$ and $q$ is a perfect square, such that
\itemitem{(1)} the Fourier transform maps each summand of $\Cal M_q^2$ onto the other,
\itemitem{(2)} the unit projection $e$ of $\Cal M_q$ is approximately central,
\itemitem{(3)} the compressions $eUe$ and $eVe$ are approximately contained in $\Cal M_q$,
\itemitem{(4)} the projection $e$ is also the unit of a subalgebra $C$
isomorphic to $M_q(C(\Bbb T))$ such that $C \oplus \sigma(C)$ is a Fourier invariant
circle algebra,
\endproclaim

\bigpagebreak

Section 2 contains the background and notation used in later sections.
In Section 3 the construction of $\Cal M_q$ is carried out to give (1) and (4).
In Sections 4 and 5 we prove (3) and (2), respectively.  Section 6 is an Appendix
containing some basic lemmas used in the approximations of Sections 4 and 5.

Throughout, we shall adopt the notation $e(t) := e^{2\pi it}$.

\bigpagebreak

\specialhead \S2.  Initial Considerations \endspecialhead

In order to focus on the C*-algebraic aspects of the problem, and to avoid some number
theoretic obstacles, we shall from now on assume that $\theta$ is a transcendental number 
in $(0,1)$.  This will simplify the choice of the lattice $D$ in Section 3.  The
following result is an application of the Thue-Siegel-Roth Theorem (see for example
[\SL], Chapter 7).

\proclaim{Proposition 2.1}
Let $\theta$ be a transcendental number in $(0,1)$.  For each $a>0$, there are 
infinitely many convergents $p/q$ such that $p$, $q$, $q-p$ are all perfect squares 
and satisfying
$$
\left| \theta - \frac pq \right| \ < \ \frac1{q^a}.
\tag2.1
$$
\endproclaim
\demo{Proof}
Since $\theta$ is transcendental, so is 
$$
\xi = \frac{1-\sqrt\theta}{\sqrt{1-\theta}}
$$
which is in $(0,1)$.  Thanks to the Thue-Siegel-Roth Theorem, for each positive number 
$\epsilon$ there are infinitely many convergents $r/s\in(0,1)$ of $\theta$ that satisfy
$$
\left| \frac rs - \xi \right|\ < \ \frac1{s^\epsilon}.
$$
Let $f(x)=2x/(x^2+1)$ so that $f(\xi)=\sqrt{1-\theta}$.  One clearly has
$$
|f(x)-f(y)| = \frac{2|x-y|\,|1-xy|}{(x^2+1)(y^2+1)}.
\tag2.2
$$
Let $m = 2rs, \ k = s^2-r^2, \ n = r^2+s^2,$
so that $m^2+k^2=n^2$ and $\tfrac mn = f(\tfrac rs)$.  Since $\xi$ and $r/s$ are in
$(0,1)$, from (2.2) we have
$$
\left| \frac mn - \sqrt{1-\theta} \right| = |f(\tfrac rs) - f(\xi)| < 
2\left|\frac rs - \xi\right| < \frac 2{s^\epsilon}.
$$
One gets
$$
\left| \theta - \frac{n^2-m^2}{n^2} \right| =
\left| \frac{m^2}{n^2} - (1-\theta) \right| < \frac4{s^\epsilon}.
$$
Now since $r<s$, we have
$$
\frac4{s^\epsilon} = \frac{4 \cdot 2^{\epsilon/2}}{(2s^2)^{\epsilon/2}} < 
\frac{4 \cdot 2^{\epsilon/2}}{(r^2+s^2)^{\epsilon/2}} = 
\frac{4 \cdot 2^{\epsilon/2}}{n^{\epsilon/2}}.
$$
Therefore, taking $p=n^2-m^2 = k^2$ and $q=n^2$ (so $q-p=m^2$) they satisfy 
$$
\left| \theta - \frac pq \right| 
< \frac{4 \cdot 2^{\epsilon/2}}{n^{\epsilon/2}}  
= \frac{4 \cdot 2^{\epsilon/2}}{q^{\epsilon/4}}  
< \frac1{q^a}
$$
where the last inequality holds, for example by taking $\epsilon = 4(1+a)$ and for large
enough $q$.  Since the greatest common divisor $d$ of $p$ and $q$ is also
a perfect square, $\tfrac pd, \tfrac qd, \tfrac {q-p}d$ are perfect squares and
satisfy the preceding inequality, so one can assume $p,q$ are relatively prime.  In
this case, it is well-known (for $q>2$) that (2.1) in itself implies that $p/q$ must be a
convergent of $\theta$ (of course, providing that $a\ge3$).  \qed
\enddemo

\bigpagebreak

If $\theta$ is algebraic, then by Roth's Theorem [\KR], for any 
$\epsilon>0$ the inequality $|\theta - \tfrac pq | < \tfrac1{q^{2+\epsilon}}$
can only be satisfied by {\it finitely} many rationals $p/q$, a result for which he 
received the 1958 Fields Medal.  Therefore, our approach here only applies to 
transcendental numbers $\theta$, and throughout the paper we shall let $p/q$ be a 
fixed convergent of $\theta$ such that $p$, $q$, and $q-p$ are perfect squares 
satisfying (2.1) with $a=3$:  
$$
\left| \theta - \frac pq \right| \ < \ \frac1{q^3}.
$$
Further, without loss of generality, we may assume that there are
infinitely many such convergents such that $\tfrac pq < \theta$ (otherwise, one replaces
$\theta$ by $1-\theta$).
(If there are infinitely many convergents $p/q$ such that $\theta < \tfrac pq$
and $p,q,q-p$ are perfect squares, then $\tfrac{q-p}q$ is also a convergent with 
$\tfrac{q-p}q < 1-\theta$ where $q-p$ is a perfect square.)

\bigpagebreak

\subhead Theta Functions \endsubhead
The main Theta functions used in this paper are
$$
\vartheta_2(z,t) = \sum_n \e^{\pi it(n+\frac 1 2)^2} \e^{i2z(n+\tfrac12)}, \qquad
\vartheta_3(z,t) = \sum_n \e^{\pi itn^2} \e^{i2zn}
$$
for $z,t\in\Bbb C$ and $\roman{Im}(t)>0$, where all summations range over the integers.
(For a classic treatment see [\WW].)

\newpage

\subhead A Gaussian Theta Function \endsubhead
In [\SWe] a Schwartz function $h$ on $\Bbb R$ was constructed that lead to the
existence of flip invariant projections $e$ in $A_\theta$ orthogonal to their 
Fourier transform (i.e., $e\sigma(e)=0$).
Since this function is important for the present paper, we shall recall briefly some of
its properties.  Let $\alpha>0, \beta>2$ be such that $\beta^2=4(\alpha^2+1)$.  Define 
$$
h(x) = e^{-\pi\alpha x^2}\ \sum_p\ e^{-\pi\alpha p^2} e^{\pi\alpha p}
e([\tfrac\beta2 p - \tfrac\beta4]x)
$$
an even real-valued Schwartz function on $\Bbb R$.  It was shown in [\SWe] that 
$$
\int_{\Bbb R} \conj{h(x)}\, \widehat h(x+\beta m) e(\beta nx)\,dx = 0 
\tag2.3
$$
for all integers $m,n$, where $\widehat h$ is the usual Fourier transform of $h$ (defined
below). It was also shown that for all real $s,t$ one has
$$
H(s,t) := \int_{\Bbb R} \conj{h(x)} h(x+s) e(xt)\,dx =
\frac1{\sqrt{2\alpha}} e(-\tfrac{st}2) e^{-\pi\alpha s^2/2} e^{-\pi t^2/(2\alpha)}
\Gamma(\tfrac t\beta, \tfrac s \beta)
\tag2.4
$$
where
$$
\Gamma(u,v) \ = \  e^{\pi\alpha/2}
\Bigl[ 
\vartheta_2(\tfrac\pi2\beta^2v, 2i\alpha)
\vartheta_3(i\tfrac\pi{2\alpha}\beta^2u, it_\alpha)
+
\vartheta_3(\tfrac\pi2\beta^2v, 2i\alpha) 
\vartheta_2(i\tfrac\pi{2\alpha}\beta^2u, it_\alpha)
\Bigr]
\tag2.5
$$
for $u,v\in\Bbb R$, and $t_\alpha = 4\alpha + \tfrac2\alpha$ ([\SWe, Lemma 4.2]).  
It is clear that $\Gamma$ is even in each variable and is real valued (in fact all of 
its terms are positive, except possibly $\vartheta_2(\tfrac\pi2\beta^2v, 2i\alpha)$).
Using the fact that
$$
2 < e^{\pi\alpha/2} \vartheta_2(0, 2i\alpha) < 2 \vartheta_3(0, 2i\alpha)
\tag2.6
$$
one has
$$
|\Gamma(u,v)| \ \le \vartheta_3(0,2i\alpha) 
\Bigl[ 2\vartheta_3(i\tfrac\pi{2\alpha}\beta^2u, it_\alpha)
	+ e^{\pi\alpha/2} \vartheta_2(i\tfrac\pi{2\alpha}\beta^2u, it_\alpha)
\Bigr].
\tag2.7
$$
We also had the period 1 functions
$$
\psi_n(t) := \sum_m\ e(\beta^2 mn/2)  e^{-\pi \beta^2 m^2/(2\alpha)} \Gamma(m,n) \,
e(mt)
\tag2.8
$$
for $n\in\Bbb Z$, which will arise in the cut down approximations below.  The 
function $\psi_0$ is positive invertible and it was shown that we have the 
inequality
$$
\frac{|\psi_n(t)|}{\psi_0(s)}  \le 1.21 \,
\frac{\vartheta_3(0, 2i\alpha)}{\vartheta_3(\tfrac\pi2, 2i\alpha)}.
$$
for all real $s,t$ and any integer $n$ ([\SWe, Lemma 4.3]). 
The right side of this inequality is a decreasing function of $\alpha$ and hence is no 
more than its value at $\alpha=\tfrac14$, hence 
$$
\frac{|\psi_n(t)|}{\psi_0(s)}   \le 1.21 \, 
\frac{\vartheta_3(0, \tfrac i2)}{\vartheta_3(\tfrac\pi2, \tfrac i2)} = 
1.21(1+\sqrt2) < 3
\tag2.9
$$
for $\alpha > \tfrac14$. Further, for $\beta > \beta_0 := 2.063$ and all real $t$ one has
$$
1 < \psi_0(t) < 5, \ \ 
\text{and \ \ $\psi_0(t) \to 2$ uniformly in $t$ as $\beta\to\infty$}.
\tag2.10
$$
([\SWe, Proposition 4.4]).

\bigpagebreak

\subhead Rieffel's Framework \endsubhead
We briefly sketch Rieffel's setup in [\MRc] and the notation to be used.  Let $M$ be 
a locally compact Abelian group and let $G=M\times \hat M$, where $\hat M$ is the dual
of $M$. Let $\frak h$ denote the Heisenberg cocycle on the group $G=M\times \hat M$ 
given by $\frak h((m,s),(m',s')) = \langle m, s' \rangle$, where 
$\langle m, s' \rangle$ is the canonical pairing $M\times\hat M \to \Bbb T$.
The Heisenberg (unitary) representation $\pi:G\to \Cal U(L^2(M))$ of $G$ is given by
$$
(\pi_{(m,s)} f)(n) \ = \ \langle n,s \rangle f(n+m)
$$
where $m,n\in M,\  s\in\hat M$, and $f\in L^2(M)$.  It has the properties
$$
\pi_x\pi_y = \frak h(x,y) \pi_{x+y} = \frak h(x,y) \conj{\frak h(y,x)} \pi_y\pi_x,
\qquad
\pi_x^* = \frak h(x,x) \pi_{-x}
$$
for $x,y\in G$.  If $D$ is a given lattice in $G$ (i.e., discrete subgroup such that
$G/D$ is compact), its covolume $|G/D|$ is the Haar measure of a fundamental domain for 
$D$ in $G$.  Its associated twisted group C*-algebra $C^*(D,\frak h)$ is the
C*-subalgebra of the bounded operators on $L^2(M)$ generated by the unitaries
$\pi_x$ for $x\in D$.  The complementary lattice of $D$ is
$$
D^\perp = \{ y\in G: \frak h(x,y) \conj{\frak h(y,x)} = 1,\ \forall x\in D\}.
$$
The C*-algebra $C^*(D^\perp,\bar\frak h)$ can be viewed as the
C*-subalgebra of bounded operators on $L^2(M)$ generated by the unitaries
$\pi_y^*$ for $y\in D^\perp$.  By Rieffel's Theorem 2.15 [\MRc], the Schwartz space 
on $M$, denoted $\Cal S(M)$, is an equivalence bimodule with
$C^*(D,\frak h)$ acting on the left and $C^*(D^\perp,\bar\frak h)$ acting on the
right by
$$
\align
af &= \int_D a(x)\pi_x(f) dx = |G/D|\,\sum_{x\in D} a(x)\pi_x(f)
\\
f b &= \int_{D^\perp} b(y) \pi_y^*(f) dy = \sum_{y\in D^\perp} b(y) \pi_y^*(f)
\endalign
$$
where $f\in\Cal S(M),\ a\in C^*(D,\frak h),\ b\in C^*(D^\perp,\bar\frak h)$,
and where the mass point measure ($dx$) on $D$ is $|G/D|$ and on $D^\perp$ it is one.
The inner products on $\Cal S(M)$ with values in the algebras $C^*(D,\frak h)$
and $C^*(D^\perp, \conj{\frak h})$ are given, respectively, by
$$
\Dinner f g = |G/D|\,  \sum_{w\in D} \Dinner f g (w) \, \pi_w,
\qquad
\Dpinner f g = \sum_{z\in D^\perp} \Dpinner f g (z) \, \pi_z^*
$$
where
$$
\align
\Dinner f g (w_1,w_2) &= \int_M f(x) \conj{g(x+w_1)}\ \conj{\inner x {w_2}{}} dx
\\
\Dpinner f g (z_1,z_2) &= \int_M \conj{f(x)} g(x+z_1) \inner x {z_2}{} dx
\endalign
$$
where $(w_1,w_2) \in D$ and $(z_1,z_2)\in D^\perp$.  These satisfy the associativity
condition
$$
\Dinner f g h = f \Dpinner g h.
$$
(See [\MRc], pages 266 and 269.)  The canonical normalized traces are given by
$$
\vartau_D\left( \sum_{w\in D} a_w\pi_w \right) = a_0,
\qquad
\vartau_{D^\perp}\left( \sum_{z\in D^\perp} b_z\pi_z^* \right) = b_0,
$$
($a_w,b_w\in \Bbb C$) and satisfy the equality
$$
\vartau_D( \Dinner f g ) = |G/D|\, \vartau_{D^\perp} (\Dpinner g f).
$$
From this it follows that if $\xi$ is a Schwartz function such that $\Dpinner\xi\xi = 1$,
then $e=\Dinner\xi\xi$ is a projection in $C^*(D,\frak h)$ of trace $|G/D|$.
In this case, one has the cut down isomorphism 
$$
\mu_\xi:\ e C^*(D,\frak h) e \to C^*(D^\perp,\bar\frak h)
$$
given by
$$
\mu_\xi(exe) = \Dpinner {\xi}{x\xi}, \qquad 
\mu_\xi^{-1}(y) = \Dinner {\xi y}\xi
$$
for $x\in C^*(D,\frak h)$ and $y\in C^*(D^\perp,\bar\frak h)$.  
Projections $e$ arising in this manner might be called {\it generalized Rieffel 
projections}, since in [\EL], Elliott and Lin showed that the classical Rieffel 
projections [\MRb] in the rotation algebra do in fact have this inner product form.

In the present paper, we shall take $M = \Bbb R \times \Bbb F$ where $\Bbb F$ is a 
finite abelian group so that $\hat M = M$ in a natural fashion.  
With this identification in mind one is permitted to define the order four 
automorphism $R:G\to G$ by $R(u;v)=(-v;u),\ u,v\in M$.  If $D$ is a lattice
subgroup of $G$ such that $R(D)=D$ (and hence $R(D^\perp)=D^\perp$), then (as
in [\FBb]) there are order four automorphisms $\sigma, \sigma'$ of
$C^*(D,\frak h)$ and $C^*(D^\perp,\bar \frak h)$, respectively, that satisfy
$$
\sigma(\pi_w) = \conj{\frak h(w,w)} \pi_{Rw}, \qquad
\sigma'(\pi_z^*) = \frak h(z,z) \pi_{Rz}^*,
$$
for $w\in D,\ z\in D^\perp$, and
$$
\sigma(\Dinner fg) = \Dinner{\hat f}{\hat g}, \qquad
\sigma'(\Dpinner fg) = \Dpinner{\hat f}{\hat g},
$$
where $\hat f$ is the usual Fourier transform of $f\in \Cal S(M)$ given by
$$
\hat f (s) = \int_M\ f(x) \conj{\inner x s {}}\ dx
$$
for $s\in M$.  Note that the isomorphism $\mu_\xi$ commutes with the Fourier transform 
on each algebra according to the commutative diagram
$$
\CD
eC^*(D)e  @> \mu_\xi >> C^*(D^\perp)    \\
@V \sigma VV           @VV \sigma' V \\
\sigma(e) C^*(D) \sigma(e)   @>  \mu_{\hat\xi}  >>    C^*(D^\perp)
\endCD
$$
which, for $x\in eC^*(D)e$, follows from
$$
\sigma'(\mu_\xi(x)) =
\sigma'( \Dpinner{\xi}{x \xi}) =
\Dpinner{\hat\xi}{\widehat{x \xi}} =
\Dpinner{\hat\xi}{\sigma(x)\hat\xi} =
\mu_{\hat\xi}(\sigma(x))
$$
as $\sigma(e) = \Dinner{\hat\xi}{\hat\xi}$.  As in [\EL], we shall write 
$``exe" := \Dpinner \xi{x\xi}$.

\newpage

\specialhead \S3. The Construction \endspecialhead

Let $M=\Bbb R \times \Bbb Z_q \times \Bbb Z_q$, and consider the
lattice $D$ in $G=M\times \hat M$ with basis
$$
D:\ \ \bmatrix \varepsilon_1 \\ \varepsilon_2 \endbmatrix 
=
\bmatrix
a & [p']_q & 0 & 0 & 0 & 0 \\
0 & 0 & 0 & a & [p']_q & 0 
\endbmatrix,
$$
where $p'$ is a positive integer such that $p=(p')^2$, where we recall that 
$\tfrac pq < \theta$, and where $a = (\theta - \tfrac pq)^{1/2}$.
This lattice is clearly invariant under $R: (w_1,w_2)\mapsto(-w_2,w_1)$.  
The algebra $C^*(D,\frak h)$ is generated by the canonical unitaries 
$U_1=\pi_{\varepsilon_1},\ U_2=\pi_{\varepsilon_2}$ whose commutation relation is
$$
U_1U_2U_1^*U_2^*
= \pi_{\varepsilon_1} \pi_{\varepsilon_2} \pi_{\varepsilon_1}^* \pi_{\varepsilon_2}^*
= \frak h \frak h^*(\varepsilon_1, \varepsilon_2)
= e(a^2+\tfrac{(p')^2}q) = e(\theta) 
$$
Therefore, $U_1,U_2$ generate the irrational rotation algebra $A_\theta$, and the Fourier
transform enjoys $\sigma(U_1)=U_2,\ \sigma(U_2)=U_1^{-1}$.  
A fundamental domain for $D$ is
$ 
[0,a)\times \Bbb Z_q \times \Bbb Z_q \times [0,a)\times \Bbb Z_q \times \Bbb Z_q
$
so its covolume is $|G/D| = a^2 q^2 = q(q\theta-p) <1$.  Put 
$$
\beta := \frac1{qa}.
$$ 
Observe that since by hypothesis $a^2< q^{-3}$, one has $\beta>\sqrt q$ so that $\beta>2$
for $q\ge4$.  The annihilator lattice $D^\perp$ can be checked to be spanned by the basis
elements
$$
D^\perp: \ \ 
\bmatrix \delta_1 \\ \delta_2 \\ \delta_3 \\ \delta_4 \endbmatrix
=
\bmatrix
\beta & [-p_0] & 0 & 0 & 0 & 0 \\
0 & 0 & 0 & \beta & [-p_0] & 0 \\
0 & 0 & [1] & 0 & 0 & 0  \\
0 & 0 & 0 & 0 & 0 & [1]  
\endbmatrix
$$
where $p_0$ is an integer such that $p_0p' = 1 \mod q$.
Let $V_j = \pi_{\delta_j}^* = \pi_{-\delta_j}$.  Then
$$
V_1V_2 = e(\mu) V_2V_1, \qquad V_3V_4 = e(\tfrac1q) V_4V_3, \qquad
V_jV_k=V_kV_j
$$
for $j=1,2$ and $k=3,4$, where $\mu = \beta^2 + \tfrac{p_0^2}q$ is in the
$\roman{GL}(2,\Bbb Z)$ orbit of $\theta$.  One has
$C^*(D^\perp,\bar\frak h) = C^*(V_j) \cong M_q(A_\mu)$.

The $D^\perp$ inner product becomes, for $f,g\in\Cal S(M)$,
$$
\align
&\Dpinner f g (n_1\delta_1 + n_2\delta_2 + n_3\delta_3 + n_4\delta_4)
\\
&\ \ \ =
\Dpinner f g \left( n_1\beta, [-n_1p_0], [n_3];\ n_2\beta, [-n_2p_0], [n_4] 
\right)
\\
&\ \ \ =
\int_{\Bbb R \times \Bbb Z_q \times \Bbb Z_q}
\!\!\!\!\conj{ f(x,[n],[m]) }  g(x+n_1\beta, [n-n_1p_0], [m+n_3])
e( xn_2\beta - \tfrac{nn_2p_0}q + \tfrac{mn_4}q ) \ dx d[n] d[m]
\\
&\ \ \ =
\frac1q \sum_{[n],[m]\in\Bbb Z_q}
e(-\tfrac{nn_2p_0}q + \tfrac{mn_4}q )
\int_{\Bbb R} \conj{f(x,[n],[m])} g(x+n_1\beta, [n-n_1p_0], [m+n_3]) e(xn_2\beta)\,dx.
\endalign
$$
Since
$$
\pi_{n_1\delta_1 + n_2\delta_2 + n_3\delta_3 + n_4\delta_4}^*
= V_1^{n_1} V_2^{n_2} V_3^{n_3} V_4^{n_4}
$$
one has (where $n_3,n_4=0,1,\dots,q-1$ and $n_1,n_2\in\Bbb Z$)
$$
\Dpinner fg = 
\sum_{n_1n_2n_3n_4} \Dpinner fg(\Sigma n_*\delta_*)
\cdot V_1^{n_1} V_2^{n_2} V_3^{n_3} V_4^{n_4}.
$$
Now suppose that $f$ and $g$ are Schwartz functions on $M$ of the form
$$
f(x,[n],[m]) = h(x) \varphi_f([n],[m]), \qquad
g(x,[n],[m]) = h'(x) \varphi_g([n],[m])
$$
where $h,h'\in\Cal S(\Bbb R)$ and $\varphi_f, \varphi_g$ are any functions.  Then
one gets
$$
\Dpinner f g (n_1\delta_1 + n_2\delta_2 + n_3\delta_3 + n_4\delta_4)
\ =\ 
\frac1q \Omega_{n_1n_2n_3n_4}(\varphi_f,\varphi_g) \cdot
[hh'](n_1\beta,n_2\beta)
$$
where
$$
\Omega_{n_1n_2n_3n_4}(\varphi_f,\varphi_g) =
\sum_{[n],[m]\in\Bbb Z_q}
e(-\tfrac{nn_2p_0}q + \tfrac{mn_4}q )
\conj{\varphi_f([n],[m])} \varphi_g([n-n_1p_0],[m+n_3])
$$
and
$$
[hh'](s,t) = \int_{\Bbb R} \conj{h(x)} h'(x+s) e(xt)\,dx.
$$
First, observe that if $h$ is chosen so that $[h\hat h](n_1,n_2) = 0$ for all
integers $n_1, n_2$, then $\Dpinner f{\hat f} = 0$.  (This follows from
$\hat f(x,[n],[m]) = \hat h(x) \, \widehat{\varphi_f}([n],[m])$.)
Consider the Schwartz function
$$
f(x,[n],[m]) = (2\alpha)^{1/4}\,h(x) \delta_q^{n-m} 
$$
where $\delta_q^k=1$ if $q|k$ and $0$ otherwise, $h$ is the Gaussian Theta function
of Section 2, and $\alpha>0$ satisfies $\beta^2=4(\alpha^2+1)$.  We have
$$
\align
\Omega_{n_1n_2n_3n_4}(\varphi_f,\varphi_f) 
&=
\sum_{m,n=0}^{q-1}
e(-\tfrac{nn_2p_0}q + \tfrac{mn_4}q )
\delta_q^{n-m} \delta_q^{n-m - n_1p_0 - n_3}
\\
&=
\delta_q^{n_1p_0 + n_3}
\sum_{n=0}^{q-1} e(-\tfrac n q (n_2p_0 - n_4) )
\\
&=
q\,\delta_q^{n_1p_0 + n_3} \delta_q^{n_2p_0 - n_4}
\endalign
$$
hence (with $n_3,n_4=0,1,\dots,q-1$ and $n_1,n_2\in\Bbb Z$)
$$
\align
\Dpinner ff 
&= \frac{\sqrt{2\alpha}}q  \, \sum_{n_1n_2n_3n_4}
q\,\delta_q^{n_1p_0 + n_3} \delta_q^{n_2p_0 - n_4}
H(n_1\beta,n_2\beta) V_1^{n_1} V_2^{n_2} V_3^{n_3} V_4^{n_4}
\\
&= 
\sqrt{2\alpha}
\sum_{n_1n_2}
H(n_1\beta,n_2\beta)\, V_1^{n_1} V_2^{n_2} V_3^{-n_1p_0} V_4^{n_2p_0}
\\
&=: b
\endalign
$$
where (since the pair $V_1,V_2$ commutes with $V_3,V_4$)
$$
b = \sqrt{2\alpha} \sum_{n_1n_2}  H(n_1\beta,n_2\beta)\, W_1^{n_1} W_2^{n_2}, 
\qquad
W_1 = V_1 V_3^{-p_0}, \qquad W_2 = V_2 V_4^{p_0}
$$
and $W_1 W_2 = e(\beta^2)\, W_2 W_1$. (Recall that $H$ is given by (2.4).)
The following lemma shows that $b$ is an invertible element (at least) for sufficiently 
large $q$.

\bigpagebreak

\proclaim{Lemma 3.1} 
One has $ \|b - \psi_0 \|\to 0 $ as $q\to\infty$ and hence $\|b-2\|\to 0$.
In particular, $b$ is invertible for sufficiently large $q$.
\endproclaim
\demo{Proof}  Since
$$
H(m\beta,n\beta) =
\frac1{\sqrt{2\alpha}} e(-\tfrac12\beta^2mn) e^{-\pi\alpha \beta^2m^2/2}
e^{-\pi\beta^2 n^2/(2\alpha)} \Gamma(n,m)
$$
we can write
$$
b = \sum_m e^{-\pi\alpha \beta^2m^2/2}\, W_1^m \rho_m
$$
where $\rho_m = \rho_m(W_2)$ can be viewed as a period one function as
$$
\rho_m(t) =
\sum_n e(-\tfrac12\beta^2mn) e^{-\pi \beta^2 n^2/(2\alpha)} \Gamma(n,m) e(nt)
= \psi_m(t-\beta^2m)
$$
where $\psi_m$ is given by (2.8).  In particular, $\rho_0=\psi_0$ is strictly positive.
From (2.9) one has $\left\| \rho_m \rho_0^{-1} \right\| < 3$ for all $m$ 
(and $\alpha > \tfrac14$), hence one obtains
$$
\| b\psi_0^{-1} - I \| \le
\sum_{m\not=0} e^{-\pi\alpha \beta^2 m^2/2} \left\| \rho_m\rho_0^{-1} \right\|
< 3 \sum_{m\not=0} e^{-\pi\alpha \beta^2 m^2/2}
= 3[ \vartheta_3(0, \tfrac i2 \alpha\beta^2) - 1] \ =: C_q
$$
and $C_q \to 0$ as $q\to\infty$.  Thus one gets
$
\| b - \psi_0 \| < C_q \|\psi_0\| < 5 C_q.
$
Since $\|\psi_0-2\| \to 0$ as $\beta\to\infty$ by (2.10), one deduces that
$\|b-2\|\to0$.  \qed
\enddemo

\bigpagebreak

Now since $\|x-1\|< \epsilon <1$ imples $\|x^{-1}-1\| < \epsilon/(1-\epsilon)$,
one gets (for large enough $q$)
$$
\| \psi_0 b^{-1} - I \| < \frac{C_q}{1-C_q}
$$
and using (2.10)
$$
\|b^{-1} - \psi_0^{-1}\| < \frac{C_q}{1-C_q} \|\psi_0^{-1}\| < \frac{C_q}{1-C_q}.
$$
This proves the following lemma to be used below.

\proclaim{Lemma 3.2}
The norms $\|b\|$, $\|b^{-1/2}\|$ are uniformly bounded for $\alpha > \tfrac14$.
\endproclaim

\bigpagebreak

Putting $\xi = fb^{-1/2}$ one has $\Dpinner \xi\xi = 1$ and hence
$
e = \Dinner \xi \xi
$
is a projection in $C^*(D,\frak h) = A_\theta$ with trace $\vartau(e) = |G/D| =
q(q\theta-p)$.  Since $b$ is flip invariant, because the flip
($V_j\mapsto V_j^{-1}$) maps $W_j$ to $W_j^{-1}$ and
$H(-m\beta, -n\beta) = H(m\beta, n\beta)$, it follows that $e$ is a flip invariant 
projection in $A_\theta$ such that $e\sigma(e)=0$ (by (2.3)).
We conclude that $e$ is the unit projection of a matrix subalgebra $\Cal M_q$
of $A_\theta$ obtained from the Morita isomorphism $\mu_\xi$ at the end of Section 2.
Hence $A_\theta$ contains the finite dimensional C*-subalgebra 
$\Cal M_q^2 = \Cal M_q \oplus \sigma(\Cal M_q)$ in which its summands are mapped onto one
another by the Fourier transform (and whose unit $e+\sigma(e)$ is Fourier invariant).
If one looks at the C*-algebra generated by $V_1,V_3,V_4$, which is a $q\times q$
matrix algebra over $C(\Bbb T)$, and transforms it by $\mu_\xi$, one obtains a circle
subalgebra of $A_\theta$ whose unit is also $e$.  This proves (1) and (4) of the Main 
Theorem.  Since from the following two sections $e$ is shown to be approximately central, 
and that the compressions $eU_1e, eU_2e$ are approximately contained in $\Cal M_q$, it 
follows that $\Cal M_q^2$ is a partially approximating C*-subalgebra of $A_\theta$ that 
is invariant under the Fourier transform.

\remark{Remark}
It is immediate from the present construction that one also obtains Fourier orthogonal 
projections of trace $k(q\theta-p)$ for $k=1,\dots,q$. (By looking at subprojections of 
$e$.) 
\endremark

\newpage

\specialhead \S4.  Cut Down Approximation \endspecialhead

In this section the compressions $eU_1e$ and $eU_2e$ are approximated by elements of
the matrix algebra generated by $V_3,V_4$ of the previous section.  More precisely:

\proclaim{Proposition 4.1} One has
$$
\| \Dpinner \xi{U_1\xi} - V_3^{p'} \| \to 0, \qquad
\| \Dpinner \xi{U_2\xi} - V_4^{-p'} \| \to 0, 
$$
as $q \to \infty$.
\endproclaim

Since $U_1= \pi_{\varepsilon_1}$, 
$$
(U_1f)(x,[r],[s]) = f(x+a,[r+p'],[s])
$$
and thus
$$
\align
&\Dpinner f{U_1f}(\Sigma n_*\delta_*)
\\
&\ \ \ \ =
\frac1q \sum_{m,n=0}^{q-1}
e(-\tfrac{nn_2p_0}q + \tfrac{mn_4}q )
\int_{\Bbb R} \conj{f(x,[n],[m])} (U_1f)(x+n_1\beta, [n-n_1p_0], [m+n_3]) 
e(xn_2\beta)\,dx
\\
&\ \ \ \ =
\frac{(2\alpha)^{\tfrac14}}q \sum_{m,n=0}^{q-1}
e(-\tfrac{nn_2p_0}q + \tfrac{mn_4}q )
\int_{\Bbb R} \conj{h(x)} \delta_q^{n-m} f(x+n_1\beta+a, [n-n_1p_0+p'], [m+n_3])    
e(xn_2\beta)\,dx
\\
&\ \ \ \ =
\frac{\sqrt{2\alpha}}q \sum_{m,n=0}^{q-1}
e(-\tfrac{nn_2p_0}q + \tfrac{mn_4}q )
\delta_q^{n-m} \, \delta_q^{n-m - n_1p_0+p' - n_3}
\int_{\Bbb R} \conj{h(x)} h(x+n_1\beta+a) e(xn_2\beta)\,dx
\\
&\ \ \ \ =
\frac{\sqrt{2\alpha}}q \delta_q^{n_1p_0-p'+n_3} \ 
\sum_{n=0}^{q-1} \ e(-\tfrac nq(n_2p_0-n_4)) \, H(n_1\beta+a, n_2\beta)
\\
&\ \ \ \ =
\sqrt{2\alpha}
\delta_q^{n_1p_0-p'+n_3} \delta_q^{n_2p_0-n_4} \, H(n_1\beta+a, n_2\beta).
\endalign
$$
Therefore,
$$
\align
\Dpinner f{U_1f}
& = 
\sqrt{2\alpha}
\sum_{n_1n_2n_3n_4} 
\delta_q^{n_1p_0-p'+n_3} \delta_q^{n_2p_0-n_4} \, 
H(n_1\beta+a, n_2\beta) V_1^{n_1} V_2^{n_2} V_3^{n_3} V_4^{n_4}
\\
&=
\sqrt{2\alpha}
\sum_{n_1n_2} 
H(n_1\beta+a, n_2\beta) V_1^{n_1} V_2^{n_2} V_3^{p'-n_1p_0} V_4^{n_2p_0}
\\
&= V_3^{p'} X 
\endalign
$$
where
$$
X = \sqrt{2\alpha} \sum_{m,n} H(m\beta+a, n\beta) W_1^m W_2^n 
$$
is a perturbation of the element $b$ of Section 3.  This gives
$$
``eU_1e" = \Dpinner {\xi}{U_1\xi} = b^{-1/2} V_3^{p'} X b^{-1/2}.
$$
One similarly calculates $\Dpinner f{U_2f}$.  Since $U_2= \pi_{\varepsilon_2}$,
$$
(U_2f)(x,r,s) = e(ax+\tfrac{rp'}q)\ f(x,r,s)
$$
we have
$$
\align
&\Dpinner f{U_2f}(\Sigma n_*\delta_*)
\\
&\ \ \ \ =
\frac1q \sum_{m,n=0}^{q-1}
e(-\tfrac{nn_2p_0}q + \tfrac{mn_4}q )
\int_{\Bbb R} \conj{f(x,[n],[m])} (U_2f)(x+n_1\beta, [n-n_1p_0], [m+n_3])
e(xn_2\beta)\,dx
\\
&\ \ \ \ =
\frac{\sqrt{2\alpha}}q e(an_1\beta)
\sum_{m,n=0}^{q-1}
e(-\tfrac{nn_2p_0}q + \tfrac{mn_4}q ) e(\tfrac{(n-n_1p_0)p'}q)
\delta_q^{n-m} \delta_q^{n-m-n_1p_0-n_3} \, H(n_1\beta, n_2\beta+a)
\\
&\ \ \ \ =
\frac{\sqrt{2\alpha}}q e(\tfrac {n_1}q) e(\tfrac{-n_1p_0p'}q)
\sum_{m,n=0}^{q-1}
e(-\tfrac{nn_2p_0}q + \tfrac{mn_4}q ) e(\tfrac {np'}q)
\delta_q^{n-m} \delta_q^{n-m-n_1p_0-n_3} \, H(n_1\beta, n_2\beta+a)
\\
&\ \ \ \ =
\frac{\sqrt{2\alpha}}q \delta_q^{n_1p_0+n_3}
\sum_{n=0}^{q-1} 
e(-\tfrac{nn_2p_0}q + \tfrac{nn_4}q + \tfrac {np'}q)\, H(n_1\beta, n_2\beta+a)
\\
&\ \ \ \ =
\frac{\sqrt{2\alpha}}q \delta_q^{n_1p_0+n_3}
\sum_{n=0}^{q-1} e(-\tfrac nq(n_2p_0-n_4-p'))\, H(n_1\beta, n_2\beta+a)
\\
&\ \ \ \ =
\sqrt{2\alpha}
\delta_q^{n_1p_0+n_3} \, \delta_q^{n_2p_0-n_4-p'}\, H(n_1\beta, n_2\beta+a)
\endalign
$$
hence
$$
\align
\Dpinner f{U_2f}
& =
\sqrt{2\alpha}
\sum_{n_1n_2n_3n_4}
\delta_q^{n_1p_0+n_3} \delta_q^{n_2p_0-n_4-p'} \,
H(n_1\beta, n_2\beta+a) V_1^{n_1} V_2^{n_2} V_3^{n_3} V_4^{n_4}
\\
&=
\sqrt{2\alpha} \sum_{n_1n_2}
H(n_1\beta, n_2\beta+a) V_1^{n_1} V_2^{n_2} V_3^{-n_1p_0} V_4^{n_2p_0-p'}
\\
&= Y V_4^{-p'} 
\endalign
$$
where
$$
Y = \sqrt{2\alpha} \sum_{mn} H(m\beta, n\beta+a) W_1^m W_2^n.
$$
Thus
$$
``eU_2e" = \Dpinner {\xi}{U_2\xi} = b^{-1/2} Y V_4^{-p'} b^{-1/2}.
$$

\bigpagebreak
\proclaim{Lemma 4.2} The norms $\|X-b\|$ and $\|Y-b\|$ converge to zero as $q\to\infty$.
\endproclaim

\demo{Proof} We shall prove this for $X$ only, the proof for $Y$ being similar.
We have
$$
\align
\|X-b\| &\le \sqrt{2\alpha} \sum_{m,n} |H(m\beta+a, n\beta) - H(m\beta, n\beta)|
\\
&= \sum_{m,n} e^{-\tfrac{\pi\beta^2n^2}{2\alpha}} e^{-\tfrac12\pi\alpha\beta^2m^2}
\, \left|e(-\tfrac12n\beta a) e^{-\tfrac12\pi\alpha a^2} e^{-\pi\alpha\beta am} 
\Gamma(n,m+\tfrac a\beta) - \Gamma(n,m)\right|
\\
&= A + 2B
\endalign
$$
where
$$
A = \sum_m e^{-\tfrac12\pi\alpha\beta^2m^2}
\, \left| e^{-\tfrac12\pi\alpha a^2} e^{-\pi\alpha\beta am}
\Gamma(0,m+\tfrac a\beta) - \Gamma(0,m)\right|
$$
and since the absolute value expression is even in $n$
$$
B = \sum_m \sum_{n=1}^\infty 
e^{-\tfrac{\pi\beta^2n^2}{2\alpha}} e^{-\tfrac12\pi\alpha\beta^2m^2}
\, \left|e(-\tfrac12n\beta a) e^{-\tfrac12\pi\alpha a^2} e^{-\pi\alpha\beta am}
\Gamma(n,m+\tfrac a\beta) - \Gamma(n,m)\right|.
$$
Using (2.7) we can estimate $B$ as follows
$$
\align
B &\le \vartheta_3(0,2i\alpha) \sum_m \sum_{n=1}^\infty
e^{-\tfrac{\pi\beta^2n^2}{2\alpha}} e^{-\tfrac12\pi\alpha\beta^2m^2}
(e^{-\tfrac12\pi\alpha a^2} e^{-\pi\alpha\beta am} + 1)
\\
&\hskip2in
\Bigl[ 2\vartheta_3(i\tfrac\pi{2\alpha}\beta^2n, it_\alpha)
       + e^{\pi\alpha/2} \vartheta_2(i\tfrac\pi{2\alpha}\beta^2n, it_\alpha)
\Bigr]
\\
&= \vartheta_3(0,2i\alpha) M\,[2N_1+e^{\pi\alpha/2}N_2]
\endalign
$$
where
$$
M = \sum_m e^{-\tfrac12\pi\alpha\beta^2m^2} 
(e^{-\tfrac12\pi\alpha a^2} e^{-\pi\alpha\beta am} + 1)
$$
and
$$
N_1 = \sum_{n=1}^\infty e^{-\tfrac{\pi\beta^2n^2}{2\alpha}} 
\vartheta_3(i\tfrac\pi{2\alpha}\beta^2n, it_\alpha),
\qquad
N_2 = \sum_{n=1}^\infty e^{-\tfrac{\pi\beta^2n^2}{2\alpha}} 
      \vartheta_2(i\tfrac\pi{2\alpha}\beta^2n, it_\alpha).
$$
We will show that $M$ is bounded over $\alpha>\tfrac14$ and that $N_1\to0$ and 
$e^{\pi\alpha/2}N_2 \to 0$ as $q\to\infty$.  For $N_1$ one uses Lemma A.2 to get
(using $t_\alpha>4$ for all $\alpha>0$)
$$
N_1 \le (1+t_\alpha^{-1/2})\, \sum_{n=1}^\infty 
e^{-\pi\tfrac{\alpha\beta^2}{2(2\alpha^2+1)}n^2} 
< \frac32 \cdot \frac12[\vartheta_3(0,\tfrac{i\alpha\beta^2}{2(2\alpha^2+1)})-1]
$$
which clearly goes to $0$ as $q\to\infty$ since
$\tfrac{\alpha\beta^2}{2(2\alpha^2+1)} \to \infty$.  For $N_2$, one writes
$$
N_2 = \sum_{n=1}^\infty e^{-\tfrac{\pi\beta^2n^2}{2\alpha}}
\sum_k e^{-\pi t_\alpha(k+\tfrac12)^2} e^{-\tfrac{\pi\beta^2n}\alpha(k+\tfrac12)}
=
\sum_k e^{-\pi t_\alpha(k+\tfrac12)^2} R(k)
$$
where, upon putting $c = \tfrac{\beta^2}{2\alpha},\ d=k+\tfrac12$, one has
$$
\align
R(k) &:= \sum_{n=1}^\infty e^{-\tfrac{\pi\beta^2n^2}{2\alpha}}
e^{-\tfrac{\pi\beta^2n}\alpha(k+\tfrac12)}
= \sum_{n=1}^\infty e^{-\pi cn^2} e^{-2\pi cdn}
= e^{\pi cd^2} \sum_{n=1}^\infty e^{-\pi c(n+d)^2}
\\
&< e^{\pi cd^2} \sum_{n=-\infty}^\infty e^{-\pi c(n+d)^2}
= e^{\pi cd^2} \sum_{n=-\infty}^\infty e^{-\pi c(n+\tfrac12)^2}
= e^{\pi cd^2} \vartheta_2(0,i\tfrac{\beta^2}{2\alpha}).
\endalign
$$
Hence
$$
N_2 <
\sum_k e^{-\pi t_\alpha(k+\tfrac12)^2} e^{\pi \tfrac{\beta^2}{2\alpha}(k+\tfrac12)^2}
\vartheta_2(0,i\tfrac{\beta^2}{2\alpha})
=
\vartheta_2(0,2i\alpha) \, \vartheta_2(0,i\tfrac{\beta^2}{2\alpha}).
$$
Therefore, by (2.6)
$$
e^{\pi\alpha/2}\,N_2
< e^{\pi\alpha/2} \vartheta_2(0,2i\alpha) \vartheta_2(0,i\tfrac{\beta^2}{2\alpha})
< 2 \vartheta_3(0,2i\alpha) \vartheta_2(0,i\tfrac{\beta^2}{2\alpha})
$$
which goes to zero (rapidly).

Next we show that $M$ is bounded.  Applying Lemma A.2,
$$
\align
M &= e^{-\tfrac12\pi\alpha a^2}
\sum_m e^{-\tfrac12\pi\alpha\beta^2m^2} e^{-\pi\alpha m/q} + 
\sum_m e^{-\tfrac12\pi\alpha\beta^2m^2}
\\
&= e^{-\tfrac12\pi\alpha a^2} 
\vartheta_3(\tfrac{i\pi\alpha}{2q}, i\tfrac12\alpha\beta^2)
+ \vartheta_3(0,i\tfrac12\alpha\beta^2)
\\
&< [1+ \tfrac{\sqrt2}{\beta\sqrt\alpha}]\, e^{-\tfrac12\pi\alpha a^2}
\exp(\tfrac{\pi^2\alpha^2}{4q^2}\cdot\tfrac2{\pi\alpha\beta^2}) 
+ \vartheta_3(0,i\tfrac12\alpha\beta^2)
\\
&= 1 + \tfrac{\sqrt2}{\beta\sqrt\alpha} + \vartheta_3(0,i\tfrac12\alpha\beta^2)
\endalign
$$
which converges to $2$ as $q\to\infty$.  Therefore, $B\to0$.  Next, we show that $A\to 0$.
We have $A = |A_1| + A_2$ where
$$
A_1 = e^{-\tfrac12\pi\alpha a^2}\Gamma(0,\tfrac a\beta) - \Gamma(0,0),
\qquad
A_2 = \sum_{m\not=0} e^{-\tfrac12\pi\alpha\beta^2m^2}
\, \left| e^{-\tfrac12\pi\alpha a^2} e^{-\pi\alpha\beta am}
\Gamma(0,m+\tfrac a\beta) - \Gamma(0,m)\right|.
$$
Using (2.6) and (2.7) again gives
$$
\align
A_2 &< \vartheta_3(0,2i\alpha) 
[ 2\vartheta_3(0, it_\alpha) + e^{\pi\alpha/2} \vartheta_2(0, it_\alpha) ]
\sum_{m\not=0} e^{-\tfrac12\pi\alpha\beta^2m^2}\, 
\left[ e^{-\pi\alpha m/q} e^{-\tfrac12\pi\alpha a^2} + 1 \right]
\\
&< 4 \vartheta_3(0,2i\alpha) \vartheta_3(0, it_\alpha)
\sum_{m\not=0} e^{-\tfrac12\pi\alpha\beta^2m^2}\, 
\left[ e^{-\pi\alpha m/q} e^{-\tfrac12\pi\alpha a^2} + 1 \right].
\endalign
$$
Expressing the summation as a $\vartheta_3$ function and using Lemma A.2 one can show
that it goes to zero as $q\to\infty$.  Lastly, for $A_1$, it will be enough to check
that
$$
\varepsilon :=
\kappa e^{\pi\alpha/2}\vartheta_2(t, 2i\alpha) - e^{\pi\alpha/2}\vartheta_2(0, 2i\alpha)
\ \to \ 0
$$
where $\kappa := e^{-\tfrac12\pi\alpha a^2} < 1$ (which converges to $1$ as $q\to\infty$) 
and $t=\tfrac\pi{2q}$.  Since
$$
e^{\pi\alpha/2}\vartheta_2(t, 2i\alpha) \ = \ e^{it} 
\sum_n e^{-2\pi\alpha n^2} e^{-2\pi\alpha n} e^{2it n}
$$
one gets
$$
\align
\varepsilon &=
\sum_n e^{-2\pi\alpha n^2} e^{-2\pi\alpha n}\, [\kappa e^{it} e^{2it n} - 1]
\\
&=
[\kappa e^{it} - 1] + [\kappa e^{-it} - 1]
+
\sum_{n\not=0,-1} e^{-2\pi\alpha n^2} e^{-2\pi\alpha n}\,[\kappa e^{it}e^{2it n}-1].
\endalign
$$
The first two terms here clearly go to zero, and the summation has absolute value no 
more than (since $n^2+2n \ge 0$ for all integers $n\not=-1$)
$$
2 \sum_{n\not=0,-1} e^{-2\pi\alpha n^2} e^{-2\pi\alpha n}
=
2 \sum_{n\not=0,-1} e^{-\pi\alpha n^2} e^{-\pi\alpha(n^2+2n)}
<
2 \sum_{n\not=0,-1} e^{-\pi\alpha n^2} \ < \ 2 [ \vartheta_3(0, i\alpha) - 1] 
$$
which clearly goes to zero.  \qed
\enddemo

\bigpagebreak

We now have
$$
\align
\| b^{-\tfrac12} V_3^{p'} X b^{-\tfrac12} - V_3^{p'} \|
&\le
\| b^{-\tfrac12} V_3^{p'} X b^{-\tfrac12} - b^{-\tfrac12} V_3^{p'} b^{\tfrac12} \|
+
\| b^{-\tfrac12} V_3^{p'} b^{\tfrac12} - b^{-\tfrac12} V_3^{p'} 2^{\tfrac12} \|
\\
& \ \ \ \ \ \ \ 
+ \| b^{-\tfrac12} V_3^{p'} 2^{\tfrac12} - 2^{-\tfrac12} V_3^{p'} 2^{\tfrac12} \|
\\
&\le
\|b^{-1}\| \|X-b\| + \| b^{-\tfrac12} \| \, \| b^{\tfrac12} - 2^{\tfrac12}\|
+ 2^{\tfrac12} \| b^{-\tfrac12} - 2^{-\tfrac12}\|.
\endalign
$$
Since the norms of $b$, $b^{-1}$, and $b^{\tfrac12}$, are uniformly bounded
by Lemma 3.2 (and Lemma A.3), and since $\|b-2\| \to 0$, by Lemma 3.1, we conclude that 
$\| \Dpinner \xi{U_1\xi} - V_3^{p'} \| \to 0$ as $q \to \infty$.  In the same fashion 
one shows that $\| \Dpinner \xi{U_2\xi} - V_4^{-p'} \| \to 0$, proving Proposition 4.1.

\newpage

\specialhead \S5.  Approximate Centrality \endspecialhead

The main result of this section is the following, which proves (2) of the Main Theorem.

\proclaim{Theorem 5.1} 
The Fourier orthogonal projection $e=\Dinner\xi\xi$ (constructed in Section 3) is 
approximately central in $A_\theta$.  That is, for $j=1,2$, one has $\|U_je-eU_j\| \to 0$
as $q\to\infty$.
\endproclaim

We begin by observing that the norm $\|U_1eU_1^*-e\|$ can be made
arbitrarily small provided that $\|U_1\Dinner ff U_1^* - \Dinner ff \|$ can be made
arbitrarily small.  This follows from the basic fact that for any vector $\eta$ in 
an equivalence $A$-$B$-bimodule and for $d,c\in B$ one has
$$
\align
\| \inner{\eta d}{\eta d} A - \inner{\eta c}{\eta c} A \|
&\le 
\| \inner{\eta d}{\eta d} A - \inner{\eta d}{\eta c} A \|  +
\| \inner{\eta d}{\eta c} A - \inner{\eta c}{\eta c} A \|
\\
&= 
\| \inner{\eta d}{\eta (d-c)} A \| +  \| \inner{\eta (d-c)}{\eta c} A \|
\\
&\le
\|\eta d \| \cdot \| \eta (d-c) \| + \|\eta (d-c) \| \cdot \| \eta c \|
\\
&\le
\|\eta\|^2\, (\|d\|+\|c\|)\, \|d-c\|
\\
&=
\|\inner \eta\eta B\| \cdot (\|d\|+\|c\|)\, \|d-c\|
\endalign
$$
where $\|\eta\| = \|\inner \eta\eta B\|^{1/2} = \|\inner \eta\eta A\|^{1/2}$ 
(by Rieffel's Proposition 3.1 of [\MRa]).  
Therefore, with $\eta=f, d=b^{-1/2}, c=2^{-1/2}$ one has
$$
\| e - \tfrac12 \Dinner f f \| \ \le\ 
\|b\| \left( \|b^{-1/2}\|+2^{-1/2} \right)\, \|b^{-1/2}-2^{-1/2}\|
$$
which goes to zero as $q\to\infty$ by Lemmas 3.1 and 3.2.
This also yields the eventual vanishing of 
$\|U_jeU_j^* - \tfrac12 U_j\Dinner f f U_j^*\|$.  Therefore, it is enough to show that 
$$
\| \Dinner ff -  U_j\Dinner f f U_j^*\| \ \to \ 0 
$$
as $q\to\infty$ for $j=1,2$, so that the approximate centrality of the projection $e$ 
follows.  

One can check that
$$
\Dinner ff = \frac{\sqrt{2\alpha}}{\beta^2} \sum_{m,n} 
\conj{H(\tfrac m\beta, \tfrac n\beta)}\ U_2^{qn} U_1^{qm}.
$$
Hence
$$
U_1\Dinner ff U_1^* - \Dinner ff = \frac{\sqrt{2\alpha}}{\beta^2} \sum_{m,n} 
\conj{H(\tfrac m\beta, \tfrac n\beta)}\ [\lambda^{qn}-1] U_2^{qn} U_1^{qm}
$$
and as $|H|$ is even in each variable one gets 
$$
\align
\|U_1\Dinner ff U_1^* - \Dinner ff\| \ 
&\le 
\frac{\sqrt{2\alpha}}{\beta^2} \sum_{m,n} 
|H(\tfrac m\beta, \tfrac n\beta)|\ |\lambda^{qn}-1|
\\
&=
\frac{2\sqrt{2\alpha}}{\beta^2} \sum_m \sum_{n=1}^\infty
|H(\tfrac m\beta, \tfrac n\beta)|\ |\lambda^{qn}-1|
\endalign
$$
Since $\lambda^{qn} = e(qn \theta) = e( qn(\theta- \tfrac pq)) = e(qna^2) = 
e(\tfrac n{q\beta^2})$ and since $|e(t)-1|\le 2\pi t$ for $t\ge0$ one has 
$|\lambda^{qn}-1| \le \tfrac{2\pi n}{q\beta^2}$.  Hence
$$
\varepsilon_1 := \|U_1\Dinner ff U_1^* - \Dinner ff\| \ \le \  
\frac{4\pi\sqrt{2\alpha}}{q\beta^4} 
\sum_m \sum_{n=1}^\infty\ n |H(\tfrac m\beta, \tfrac n\beta)|.
$$
Inserting
$$
H(\tfrac m\beta, \tfrac n\beta) = 
\frac1{\sqrt{2\alpha}} e(-\tfrac{mn}{2\beta^2}) e^{-\tfrac{\pi\alpha}{2\beta^2}m^2} 
e^{-\tfrac\pi{2\alpha\beta^2} n^2} \Gamma(\tfrac n{\beta^2}, \tfrac m{\beta^2})
$$
and using the inequality (from (2.7))
$$
|\Gamma(\tfrac n{\beta^2}, \tfrac m{\beta^2})| \ \le \ 
\vartheta_3(0,2i\alpha) \Bigl[ 2\vartheta_3(\tfrac{i\pi n}{2\alpha}, it_\alpha) +
e^{\pi\alpha/2} \vartheta_2(\tfrac{i\pi n}{2\alpha}, it_\alpha) \Bigr]
$$
yields
$$
\align
\varepsilon_1
&\le
\frac{4\pi\vartheta_3(0,2i\alpha)}{q\beta^4}
\left[\sum_m e^{-\tfrac{\pi\alpha}{2\beta^2}m^2} \right] 
\sum_{n=1}^\infty\ n
e^{-\tfrac\pi{2\alpha\beta^2} n^2} 
\Bigl[ 2\vartheta_3(\tfrac{i\pi n}{2\alpha}, it_\alpha) +
e^{\pi\alpha/2} \vartheta_2(\tfrac{i\pi n}{2\alpha}, it_\alpha) \Bigr]
\\
&=
\frac{4\pi}{q\beta^4}
\vartheta_3(0,2i\alpha)\vartheta_3(0, i\tfrac\alpha{2\beta^2}) 
\sum_{n=1}^\infty\ n
e^{-\tfrac\pi{2\alpha\beta^2} n^2} 
\Bigl[ 2\vartheta_3(\tfrac{i\pi n}{2\alpha}, it_\alpha) +
e^{\pi\alpha/2} \vartheta_2(\tfrac{i\pi n}{2\alpha}, it_\alpha) \Bigr]
\\
&=
\frac{4\pi}{q\beta^4}
\vartheta_3(0,2i\alpha)\vartheta_3(0, i\tfrac\alpha{2\beta^2})
(2A_1 + B_1)
\endalign
$$
where
$$
A_1 = \sum_{n=1}^\infty\ n
e^{-\tfrac\pi{2\alpha\beta^2} n^2} \vartheta_3(\tfrac{i\pi n}{2\alpha}, it_\alpha),
\qquad
B_1 = e^{\pi\alpha/2} \sum_{n=1}^\infty\ n
e^{-\tfrac\pi{2\alpha\beta^2} n^2} \vartheta_2(\tfrac{i\pi n}{2\alpha}, it_\alpha).
$$
These are estimated by the following lemma, in the proof of which we make repeated 
implicit use of Lemmas A.1 and A.2 of the Appendix.

\proclaim{Lemma 5.2} One has
$$
A_1 < \frac{9\beta^2}{2\pi\alpha}(2\alpha^2+1), \qquad
B_1 < \frac{12}\pi \alpha\beta^2 \vartheta_3(0,2i\alpha) +
6\beta^2(2+\beta\sqrt{2\alpha}).
$$
\endproclaim
\demo{Proof} For $A_1$, using Lemma A.2 one has
$$
\vartheta_3(\tfrac{i\pi n}{2\alpha}, it_\alpha) \ < \ 
\frac32 \exp( \tfrac{\pi n^2}{4\alpha^2 t_\alpha} )
$$
thus
$$
A_1 < \frac32 \sum_{n=1}^\infty n e^{-\pi cn^2} < 
\frac32 \left( \frac1{2\pi c} + \frac2{\sqrt{2\pi ec}} \right)
$$
where $c = -\tfrac1{2\alpha\beta^2} + \tfrac1{4\alpha^2 t_\alpha} = \tfrac \alpha
{2\beta^2(2\alpha^2+1)}$.  Since $2\pi c<1$ (as $c<\tfrac18$) and one obtains
$\tfrac1{\sqrt{2\pi c}} < \tfrac1{2\pi c}$ one gets $A_1 < \tfrac9{4\pi c} = 
\frac{9\beta^2}{2\pi\alpha}(2\alpha^2+1)$.

To estimate $B_1$ a little more work is needed.  One writes
$$
B_1 = e^{\pi\alpha/2} \sum_{n=1}^\infty\ n
e^{-\tfrac{\pi n^2}{2\alpha\beta^2}} 
\sum_k e^{-\pi t_\alpha(k+\tfrac12)^2} e^{-\tfrac{\pi n}\alpha(k+\tfrac12)}
=
e^{\pi\alpha/2} \sum_k e^{-\pi t_\alpha(k+\tfrac12)^2} N(k)
$$
where
$$
N(k) = \sum_{n=1}^\infty\ n e^{-\tfrac{\pi n^2}{2\alpha\beta^2}} 
e^{-\tfrac{\pi n}\alpha(k+\tfrac12)}.
$$
To estimate $N(k)$, put $a=\tfrac1{2\alpha\beta^2}$ and $b=\beta^2(k+\tfrac12)$ to get
$$
\align
N(k) &= e^{\pi ab^2} \sum_{n=1}^\infty n e^{-\pi a(n+b)^2}
\\
&< e^{\pi ab^2} \left[ 2|b| + \sqrt{\frac2{\pi a}} + \int_0^\infty xe^{-\pi a(x+b)^2}\,
dx \right]
\\
&< e^{\pi ab^2} \left[ 2|b| + \sqrt{\frac2{\pi a}} + \frac1{\pi a}+ \frac{|b|}{\sqrt a}
\right].
\endalign
$$
As $\alpha>\tfrac14$ and $\beta>2$, $a<\tfrac12$, and thus $\sqrt{\frac2{\pi a}} <
\frac2{\pi a}$ and we have
$$
N(k) < e^{\tfrac{\pi\beta^2}{2\alpha}(k+\tfrac12)^2} \ \left[ \frac3{\pi a} + |b|
\left(2+\frac1{\sqrt a}\right) \right].
$$
Since $-t_\alpha + \tfrac{\beta^2}{2\alpha} = -2\alpha$, one gets 
$$
\align
B_1 &< e^{\pi\alpha/2} \sum_k e^{-2\pi\alpha(k+\tfrac12)^2} 
\left[ \frac6\pi \alpha\beta^2 + \beta^2 (2+\beta\sqrt{2\alpha}) |k+\tfrac12| \right]
\\
&= 
\frac6\pi \alpha\beta^2 e^{\pi\alpha/2} \vartheta_2(0,2i\alpha) + 
\beta^2 (2+\beta\sqrt{2\alpha}) F(\alpha)
\endalign
$$
where
$$
\align
F(\alpha) &= e^{\pi\alpha/2} \sum_k |k+\tfrac12| e^{-2\pi\alpha(k+\tfrac12)^2}
= 2 e^{\pi\alpha/2} \sum_{k=0}^\infty (k+\tfrac12) e^{-2\pi\alpha(k+\tfrac12)^2}
\\
&= 2 \sum_{k=0}^\infty (k+\tfrac12) e^{-2\pi\alpha(k^2+k)}
\\
&< 2 \sum_{k=0}^\infty (k+\tfrac12) e^{-2\pi\alpha k^2}
= 1 + 2 \sum_{k=1}^\infty (k+\tfrac12) e^{-2\pi\alpha k^2}
\\
&< 1 + 4 \sum_{k=1}^\infty k e^{-2\pi\alpha k^2}
< 1 + 4\left[ \frac1{4\pi\alpha} + \frac1{\sqrt{\pi e\alpha}} \right]
< 6
\endalign
$$
since $\alpha>\tfrac14$.  Together with the inequality 
$e^{\pi\alpha/2} \vartheta_2(0,2i\alpha) < 2 \vartheta_3(0,2i\alpha)$
one obtains the desired estimate on $B_1$.  \qed
\enddemo

In view of this lemma, one now obtains
$$
\varepsilon_1 
< \frac{4\pi}{q\beta^2}
\vartheta_3(0,2i\alpha)\vartheta_3(0, i\tfrac\alpha{2\beta^2})
\left[
\frac9{2\pi\alpha}(2\alpha^2+1) + 
\frac{12}\pi \alpha \vartheta_3(0,2i\alpha) + 6(2+\beta\sqrt{2\alpha})
\right]
$$
and using the inequality $\vartheta_3(0, i\tfrac\alpha{2\beta^2}) \le 1 + 
\tfrac{\beta\sqrt2}{\sqrt\alpha}$, it is now evident that all terms involving $\alpha$
and $\beta$ are bounded over the interval $\alpha>\tfrac14$, and therefore this upper
bound vanishes as $q\to\infty$.  To complete the proof of Theorem 5.1, we show (by
similar ideas, but which we include for completeness) that 
$\varepsilon_2 = \|U_2 \Dinner ff U_2^* - \Dinner ff \|$ also vanishes as $q\to\infty$.
We have
$$
\align
\varepsilon_2 &\le 
\frac{\sqrt{2\alpha}}{\beta^2} \sum_{m,n}
|H(\tfrac m\beta, \tfrac n\beta)|\ |\lambda^{-qm}-1|
=
\frac{2\sqrt{2\alpha}}{\beta^2} \sum_n \sum_{m=1}^\infty
|H(\tfrac m\beta, \tfrac n\beta)|\ |\lambda^{qm}-1|
\\
&=
\frac2{\beta^2} \sum_n \sum_{m=1}^\infty
e^{-\tfrac{\pi\alpha}{2\beta^2}m^2} e^{-\tfrac\pi{2\alpha\beta^2} n^2} 
|\Gamma(\tfrac n{\beta^2}, \tfrac m{\beta^2})|\ |\lambda^{qm}-1|
\\
&\le
\frac{4\pi}{q\beta^4}\ \sum_n \sum_{m=1}^\infty m
e^{-\tfrac{\pi\alpha}{2\beta^2}m^2} e^{-\tfrac\pi{2\alpha\beta^2} n^2} 
|\Gamma(\tfrac n{\beta^2}, \tfrac m{\beta^2})|
\\
&\le
\frac{4\pi}{q\beta^4} \vartheta_3(0,2i\alpha) 
\left( \sum_{m=1}^\infty m e^{-\tfrac{\pi\alpha}{2\beta^2}m^2} \right) 
\sum_n e^{-\tfrac\pi{2\alpha\beta^2} n^2}
\Bigl[ 2\vartheta_3(\tfrac{i\pi n}{2\alpha}, it_\alpha) +
e^{\pi\alpha/2} \vartheta_2(\tfrac{i\pi n}{2\alpha}, it_\alpha) \Bigr]
\endalign
$$
As before, let
$$
A_2 = \sum_n e^{-\tfrac\pi{2\alpha\beta^2} n^2} \vartheta_3(\tfrac{i\pi n}{2\alpha},
it_\alpha), 	\qquad
B_2 = e^{\pi\alpha/2} \sum_n e^{-\tfrac\pi{2\alpha\beta^2} n^2} 
\vartheta_2(\tfrac{i\pi n}{2\alpha}, it_\alpha).
$$
First note that by Lemma A.1,
$$
\frac\alpha{\beta^2} \sum_{m=1}^\infty m e^{-\tfrac{\pi\alpha}{2\beta^2}m^2} < 
\frac\alpha{\beta^2} \left[\frac{\beta^2}{\pi\alpha} + \frac{2\beta}{\sqrt{\pi e\alpha}}
\right] < \frac1\pi + \frac{\sqrt\alpha}\beta < 2
$$
for $\alpha>\tfrac14$.  Thus, $\varepsilon_2\to0$ follows once it is shown that
$\frac1{\alpha\beta^2}(2A_2+B_2)$ is bounded for $\alpha>\tfrac14$
(and $\vartheta_3(0,2i\alpha)$ being already bounded over this range).
Since $t_\alpha = 4\alpha + \tfrac2\alpha > 4$ for $\alpha>0$,
$$
A_2 = \vartheta_3(0, it_\alpha) + 2 \sum_{n=1}^\infty 
e^{-\tfrac\pi{2\alpha\beta^2} n^2} \vartheta_3(\tfrac{i\pi n}{2\alpha}, it_\alpha)
< \vartheta_3(0, 4i) + 2A_1 
< \vartheta_3(0, 4i) + \frac{9\beta^2}{\pi\alpha}(2\alpha^2+1)
$$
(by Lemma 5.2), thus clearly $\frac{A_2}{\alpha\beta^2}$
is bounded.  Finally, 
$$
B_2 = e^{\pi\alpha/2} \vartheta_2(0, it_\alpha) +
2 e^{\pi\alpha/2} \sum_{n=1}^\infty e^{-\tfrac\pi{2\alpha\beta^2} n^2} 
\vartheta_2(\tfrac{i\pi n}{2\alpha}, it_\alpha)
\tag*
$$
and upon inserting the summation for $\vartheta_2$ as we have done above for $N(k)$
(which this time does not include the '$n$' in front of the exponential), one gets that 
the summation in (*) is approximately no more than 
$\beta\sqrt{2\alpha} \vartheta_2(0,2i\alpha)$.
Hence the second term in (*) is no more than
$2 \beta\sqrt{2\alpha} e^{\pi\alpha/2} \vartheta_2(0,2i\alpha) < 
4 \beta\sqrt{2\alpha} \vartheta_3(0,2i\alpha)$, 
which now makes it clear that $\frac{B_2}{\alpha\beta^2}$ is bounded (in fact goes to
zero) over $\alpha>\tfrac14$.  This proves that $\varepsilon_2\to0$ as $q\to\infty$ and
completes the proof of Theorem 5.1.

\newpage

\specialhead \S6. Appendix \endspecialhead

The following lemmas are easy but are included for completeness.

\proclaim{Lemma A.1} For all $a>0$ and $b$ real, one has
$$
\align
\sum_{n=1}^\infty n e^{-\pi a(n+b)^2} \ &< \ 
\sqrt{\frac2{\pi a}} + \frac1{\pi a} + 2|b| + \frac{|b|}{\sqrt a}
\\
\sum_{n=1}^\infty n e^{-\pi an^2} \ &< \ \frac1{2\pi a} + \frac2{\sqrt{2\pi ea}}.
\endalign
$$
\endproclaim
\demo{Proof}
Let $f(x) = xe^{-\pi a(x+b)^2}$.  Then $f$ is increasing on $[0,x_0]$ and decreasing on
$[x_0,\infty)$, where $x_0=\tfrac12(\sqrt{b^2+\tfrac2{\pi a}} - b)$.  Let $N\ge0$ be the
greatest integer $\le x_0$.  Then
$$
\sum_{n=1}^{N-1} f(n) < \int_0^N f(x)dx, 	\qquad
\sum_{n=N+2}^\infty f(n) < \int_{N+1}^\infty f(x)dx.
$$
Hence
$$
\sum_{n=1}^\infty f(n) < \int_0^N f(x)dx + \int_{N+1}^\infty f(x)dx + f(N) + f(N+1)
< 2f(x_0) + \int_0^\infty f(x)dx.  
$$
The first inequality holds from $2f(x_0) < 2x_0 < 2|b| + \sqrt{\tfrac2{\pi a}}$ and 
$$
\align
\int_0^\infty xe^{-\pi a(x+b)^2}dx &< \int_{-\infty}^\infty xe^{-\pi a(x+b)^2}dx 
= \int_{-\infty}^\infty |x-b|e^{-\pi ax^2}dx 
\\
&\le \int_{-\infty}^\infty |x|e^{-\pi ax^2}dx + 
|b|\,\int_{-\infty}^\infty e^{-\pi ax^2}dx
\\
&= 2\int_0^\infty xe^{-\pi ax^2}dx + \frac{|b|}{\sqrt a}
= \frac1{\pi a} + \frac{|b|}{\sqrt a}. 
\endalign
$$
Similarly, for the second equality note that the function $xe^{-\pi ax^2}$ has the
maximum value $\frac1{\sqrt{2\pi ea}}$ so that
$$
\sum_{n=1}^\infty n e^{-\pi an^2} < 
\frac2{\sqrt{2\pi ea}} + \int_0^\infty xe^{-\pi ax^2}dx 
= \frac2{\sqrt{2\pi ea}} + \frac1{2\pi a}.
\eqno\qed
$$ 
\enddemo

\newpage

\proclaim{Lemma A.2}
For $A>0$ and any real number $B$ one has
$$
0 < \vartheta(iB,iA) \le \left( 1 + \frac1{\sqrt A}\right) e^{B^2/(\pi A)}
$$
for $\vartheta=\vartheta_2,\vartheta_3$.
\endproclaim
\demo{Proof} Let $C=B/(\pi A)$, let $N$ be the greatest integer $\le C$, and put
$$
\xi = \min(C-N,\,N+1-C) \le \tfrac12.
$$
From a graphical consideration of areas one has
$$
\sum_p e^{-\pi A(p-C)^2}\ \le\ e^{-\pi A\xi^2} + \int_{\Bbb R}\ e^{-\pi A(x-C)^2} dx
\ \le\  1 + \frac1{\sqrt A}
$$
which holds for any real $C$.  Hence,
$$
\vartheta_3(iB,iA) = \sum_p e^{-\pi Ap^2} e^{2Bp} =
e^{B^2/(\pi A)} \sum_p e^{-\pi A(p-C)^2}
\le \left( 1 + \frac1{\sqrt A}\right) e^{B^2/(\pi A)}
$$
A similar argument holds for $\vartheta_2$, hence the result.  \qed
\enddemo

\bigpagebreak

\proclaim{Lemma A.3}
If $x$ and $y$ are positive invertible elements in a unital C*-algebra, then
$$
\| x^{\tfrac12} - y^{\tfrac12} \| \ \le \sqrt2 m^{-2} M^{\tfrac32} \|x-y\|
$$
where $m=\min(\|x^{-1}\|^{-1}, \|y^{-1}\|^{-1})$ and $M=\max(\|x\|, \|y\|)$.
\endproclaim
\demo{Proof}
Since the smallest and largest points of the spectrum of $x$ are $\|x^{-1}\|^{-1}$
and $\|x\|$, respectively, the spectra of $x$ and $y$ are contained in the
interval $[m,M]$.  Let $\epsilon, c$ be arbitrary positive numbers with
$\epsilon < m$ and $M<c$.  If $\gamma$ is the circle centered at $c$ of radius
$c-\epsilon$, then it is contained in the right half of the complex plane and
surrounds $[m,M]$.  Taking the principal branch of the square root function, one has
$$
x^{\tfrac12} - y^{\tfrac12} =
\frac1{2\pi i} \int_\gamma z^{\tfrac12} [(z-x)^{-1} - (z-y)^{-1}] dz
=
\frac1{2\pi i} \int_\gamma z^{\tfrac12} (z-x)^{-1} (x-y) (z-y)^{-1} dz
$$
and since $|z^{\tfrac12}|\le (2c-\epsilon)^{\tfrac12}$ and
$\|(z-x)^{-1}\| \le (m-\epsilon)^{-1}$ (and also for $y$) for all $z$ on $\gamma$
and $\int_\gamma |dz| = 2\pi(c-\epsilon)$, one obtains
$$
\| x^{\tfrac12} - y^{\tfrac12} \| \le
(m-\epsilon)^{-2} (c-\epsilon) (2c-\epsilon)^{\tfrac12} \|x-y\|.
$$
Letting $\epsilon\to0$ and $c\to M$ yields the result.  \qed
\enddemo

\bigpagebreak


\enddocument